\documentclass[10pt]{article}
\usepackage{graphicx,epsfig}
\usepackage{amsmath}

 \usepackage{color}      

\usepackage{url}

\usepackage{tabu}

\usepackage{array}
\usepackage{makecell}

\usepackage{multirow}

\usepackage{subcaption}
\usepackage[driverfallback=dvipdfm,
bookmarks=true,    
bookmarksopen=true,
bookmarksopenlevel=\maxdimen,
bookmarksdepth=10
]{hyperref}
\hypersetup{
pdfcenterwindow=true, 
pdfpagelabels=true,
pagebackref=true,            
hypertexnames=true,
plainpages=false,
unicode=false,     
pdftoolbar=true,                
    pdfmenubar=true,          
    pdffitwindow=false,         
    pdfstartview={FitH},        
    pdftitle={A Flux-Correction Form of the Third-Order Edge-Based Scheme for a General Numerical Flux Function},
    pdfauthor={Hiroaki Nishikawa},    
    pdfsubject={A Flux-Correction Form of the Third-Order Edge-Based Scheme for a General Numerical Flux Function},       
    pdfcreator={Hiroaki Nishikawa},   
    pdfproducer={Hiroaki Nishikawa}, 
    baseurl={https://www.researchgate.net/profile/Hiroaki-Nishikawa-2},
    pdfkeywords={}, 
    pdfnewwindow=true,     
    colorlinks=true,       
    linkcolor=black,       
    citecolor=black,       
    filecolor=black,       
    urlcolor=black,        
  breaklinks=true,
  hyperfigures=true,
  backref=true,
breaklinks=false, 
pdfpagelabels,      
pagebackref,        
hypertexnames=true, 
plainpages=false,   
naturalnames        
}
\usepackage[all]{hypcap}

\numberwithin{equation}{section}
\numberwithin{subsubsection}{subsection}
\numberwithin{subsection}{section}

\usepackage[thinlines]{easytable}

\usepackage{cleveref}
\usepackage[thinlines]{easytable}

\usepackage{blindtext,graphicx}
\usepackage[absolute]{textpos}
\numberwithin{equation}{section}
\numberwithin{subsubsection}{subsection}
\numberwithin{subsection}{section}

\definecolor{black}{rgb}{0.43, 0.21, 0.1} 
\usepackage{enumitem,xcolor}

\begin{document}


\title{\bf A Flux-Correction Form of the Third-Order Edge-Based Scheme for a General Numerical Flux Function}
 
 \author{
{Hiroaki Nishikawa}, 
hiro@nianet.org \\
  {\itshape {National Institute of Aerospace}, Hampton, VA 23666, USA}
}

\date{} 

\maketitle 

\vspace{-0.95cm}

\begin{abstract}
In this short note, we present a flux-correction form of the third-order edge-based scheme for the Euler equations that enables the direct use of a general flux function. The core idea is to replace, without loss of accuracy, the arithmetic average of the flux extrapolations by a general numerical flux evaluated at the edge midpoint, together with a correction term. We show that the proposed flux-correction form preserves third-order accuracy, provided that the general numerical flux is evaluated with the left and right states that are computed exactly for a quadratic function, which can be achieved effectively by the U-MUSCL scheme with $\kappa=1/2$. Numerical results are presented to verify third-order accuracy with the HLLC and LDFSS flux functions on irregular tetrahedral grids. 
\end{abstract}

\vspace{-0.5cm}
\section{Introduction}
\label{intro}
\vspace{-0.25cm}

\indent

The third-order edge-based scheme, originally discovered in Refs.~\cite{Katz_Sankaran_JCP:2011,diskin_thomas:AIAA2012-0609}, is an economical high-order discretization that achieves third-order accuracy on arbitrary tetrahedral grids without curved high-order grids \cite{liu_nishikawa_aiaa2016-3969,nishikawa_boundary_formula:JCP2015} or second derivatives \cite{nishikawa_liu_source_quadrature:jcp2017}. Because of these special properties, it has been the subject of active CFD algorithm research toward automated CFD simulations with anisotropic viscous grid adaptation that can be performed efficiently with simplex-element grids \cite{Kleb_etal_aiaa2019-2948,ThompsonNishikawaPadway_aiaa_scitech2023,MoriscoNishikawa_aiaa_scitech2025-0302}. The method has been formulated with a particular form of an upwind flux, which is the sum of the arithmetic average of linearly extrapolated fluxes and a dissipation term. In principle, it is possible to employ any numerical flux function by rewriting it in the form of an averaged flux and a dissipation term and applying the linearly extrapolated fluxes (see, e.g., Ref.~\cite{Fleischmann_LowDissCarbuncle_HLLC:JCP2020} for the HLLC flux). However, it will require a significant effort if the flux function of interest incorporates modifications and tuning parameters, having gone through years of development, especially for complex flow applications such as chemically-reacting hypersonic flows. In order to enable the direct use of a general flux function in the third-order edge-based scheme, we propose a flux-correction form, where any numerical flux function can be employed, and third-order accuracy is achieved by adding a correction term. The objective of this short note is to present a derivation of a flux-correction form of the third-order edge-based scheme and to demonstrate that third-order accuracy is preserved on irregular tetrahedral grids with the HLLC and LDFSS flux functions, as examples. The proposed technique is demonstrated for the Euler equations but it is directly applicable to general hyperbolic conservation laws.

\vspace{-0.35cm}
\section{Third-Order Edge-Based Scheme}
\label{intro}
\vspace{-0.25cm}

\indent

\begin{figure}[t]
\begin{center}
\begin{minipage}[b]{0.8\textwidth}
\begin{center}
        \includegraphics[width=0.39\textwidth,trim=0 0 0 0,clip]{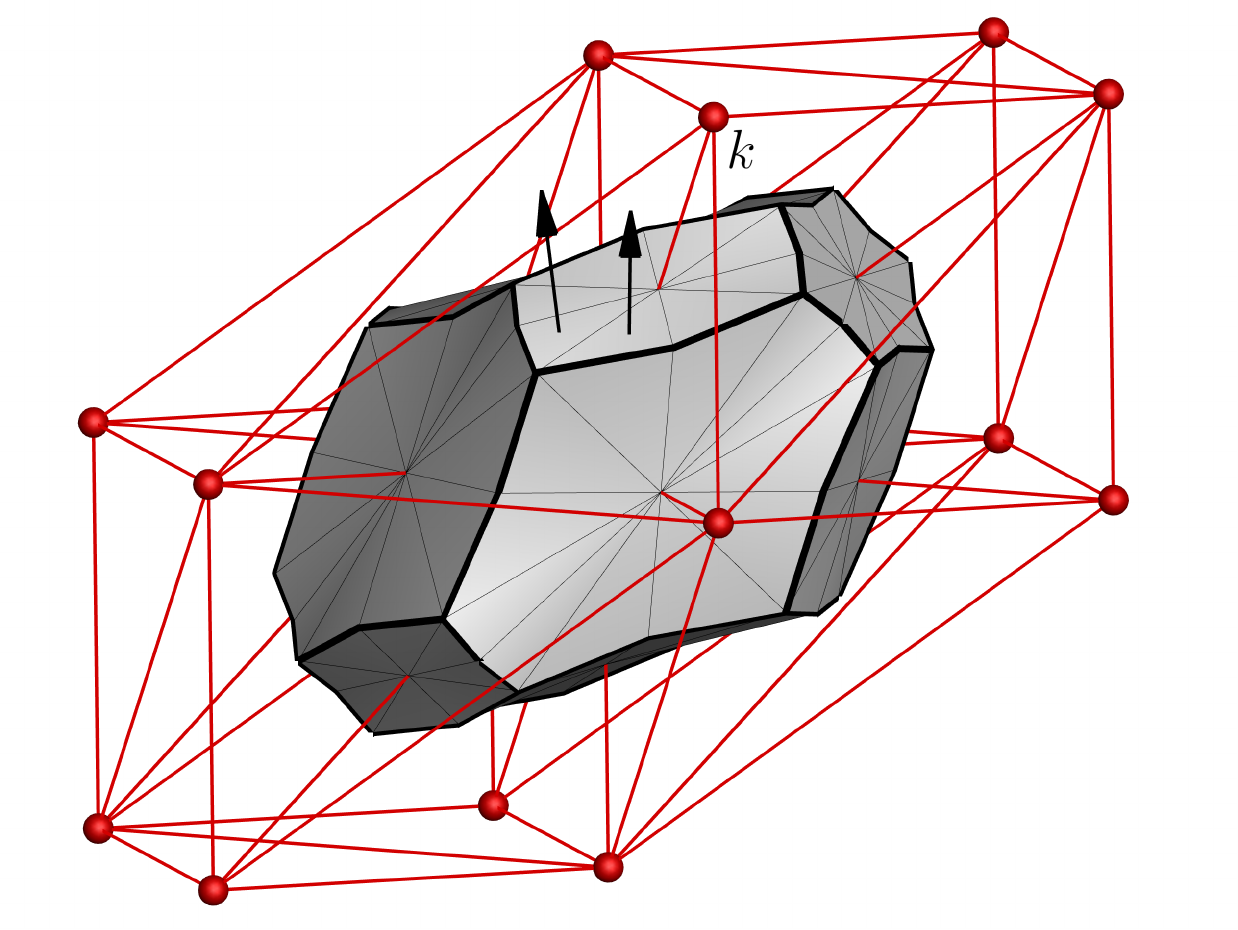}
          \caption{Tetrahedral grid and the median dual volume centered at a node $j$.}
          \label{fig:dual_vol_3d}
\end{center}
\end{minipage}
\end{center}
\end{figure}

We consider the steady Euler equations:
\vspace{-0.25cm}
\begin{eqnarray}
 \mbox{div}{\cal F}({\bf u}) = {\bf s}, 
\label{diff_form}
\end{eqnarray}
where ${\bf u}$ is a vector of conservative variables, ${\cal F}({\bf u})$ is a flux tensor, and ${\bf s}$ is a source (or forcing) term vector. In the edge-based scheme, we discretize the steady Euler equations at a node in a conservative manner over the median dual control volume defined by connecting edge midpoints and geometric centroids of elements and geometric centroids of element faces, with point-valued solutions stored at nodes. For a tetrahedral grid illustrated in Figure \ref{fig:dual_vol_3d}, the third-order edge-based scheme at a node $j$ is given by 
\vspace{-0.15cm}
\begin{eqnarray}
  \sum_{k \in \{ k_j\} } {\Phi}_{jk} | {\bf n}_{jk} |
    = 
   \frac{1}{4 D (D+2) }  \sum_{k \in \{ k_j\} } \left[  (3D+4)    {\bf s}_j  + D (  \nabla  {\bf s}_j \cdot  \Delta {\bf x}_{jk} ) -  D {\bf s}_k \right]  \left( \Delta {\bf x}_{jk} \cdot {\bf n}_{jk}  \right),
\label{threed_fv_semidiscrete_system_00}
\end{eqnarray}
\vspace{-0.35cm}
\newline
where $\{ k_j\}$ is a set of edge-connected neighbor nodes of the node $j$, ${\Phi}_{jk}$ is a numerical flux, $\Delta {\bf x}_{jk} = {\bf x}_k - {\bf x}_j $, and ${\bf n}_{jk}$ is a lumped directed-area vector, which is the sum of the directed-area vectors of the median dual faces associated with all elements sharing the edge $\{ j, k \}$ as illustrated in Figure \ref{fig:dual_vol_3d} (arrows are only included for two of the eight relevant median dual triangular faces). The right-hand side is an accuracy-preserving source quadrature formula required for third-order accuracy, where $D=1$, 2, and 3, in one, two, and three dimensions, respectively, which is one of the most efficient formulas called the compact formula derived in Ref.~\cite{nishikawa_liu_source_quadrature:jcp2017}. In this paper, we consider three dimensions and set $D=3$. The numerical flux is computed at the edge midpoint, typically by an upwind flux of the form,
\vspace{-0.4cm}
\begin{eqnarray}
 {\Phi}_{jk} ( {\bf u}_L, {\bf u}_R,  \hat{\bf n}_{jk} ) = 
  \frac{1}{2} \left[     {\bf f}_L + {\bf f}_R   \right]   
    - \frac{1}{2} {\bf D_{s}}  \left(  {\bf u}_R - {\bf u}_L \right ) ,  
  \label{upwind_flux}
\end{eqnarray}
\vspace{-0.35cm}
\newline
where ${\bf D_{s}}$ is a dissipation matrix, ${\bf f} = {\cal F} \cdot  \hat{\bf n}_{jk}$, $ \hat{\bf n}_{jk} =  {\bf n}_{jk}  /  | {\bf n}_{jk} | $, and the subscripts $L$ and $R$ indicate values extrapolated from $j$ and $k$, respectively, to the edge midpoint. The left and right states, ${\bf u}_L$ and ${\bf u}_R$, are computed from the primitive variables ${\bf w} = (\rho , u, v, w, p)$, where $\rho$ is the density, $(u,v,w)$ is the velocity vector, $p$ is the pressure: ${\bf u}_L = {\bf u}( {\bf w}_L)$ and ${\bf u}_R = {\bf u}( {\bf w}_R)$. Here, the left and right primitive variables are computed by the U-MUSCL scheme \cite{burg_umuscl:AIAA2005-4999}: 
\vspace{-0.2cm}
\begin{eqnarray}
   {\bf w}_L  =  \kappa \frac{  {\bf w}_j + {\bf w}_k}{2} + (1-\kappa) \left( {\bf w}_j   + \frac{1}{2} \nabla {\bf w}_j  \cdot \Delta {\bf x}_{jk} \right) , \quad
   {\bf w}_R =   \kappa \frac{ {\bf w}_j + {\bf w}_k}{2} + (1-\kappa) \left(   {\bf w}_k   - \frac{1}{2} \nabla {\bf w}_k \cdot \Delta {\bf x}_{jk} \right) , 
   \label{umuscl_w}
\end{eqnarray}
\vspace{-0.35cm}
\newline
where $\kappa$ is an arbitrary constant, and the gradients, $\nabla {\bf w}_j$ and $\nabla {\bf w}_k$, are computed by a quadratic least-squares fit \cite{nishikawa_onetwothree_diffusion:JCP2014} or by the implicit edge-based gradient method \cite{Nishikawa_aiaa_scitech2024}. As shown in Ref.~\cite{PadwayNishikawa:AIAAJ2022}, the jump ${\bf u}_R - {\bf u}_L$ vanishes for a quadratic function on an arbitrary grid, provided that the gradients are computed exactly for a quadratic function. Therefore, the value of $\kappa$ is arbitrary and can be adjusted to control the level of dissipation (e.g., $\kappa=1$ corresponds to zero dissipation) while preserving third-order accuracy. On the other hand, the left and right fluxes, ${\bf f}_L$ and ${\bf f}_R$, must be computed by
\vspace{-0.25cm}
\begin{eqnarray}
   {\bf f}_L  = {\bf f}_j + \frac{1}{2} \left( \frac{\partial {\bf f} }{\partial {\bf w}} \right)_{\!\! j} \nabla {\bf w}_j  \cdot \Delta {\bf x}_{jk} , \quad
   {\bf f}_R = {\bf f}_k - \frac{1}{2} \left( \frac{\partial {\bf f} }{\partial {\bf w}} \right)_{\!\! k} \nabla {\bf w}_k \cdot \Delta {\bf x}_{jk} .
\end{eqnarray}
\vspace{-0.35cm}
\newline
As proved in Ref.~\cite{nishikawa_liu_source_quadrature:jcp2017}, the edge-based scheme achieves third-order accuracy by generating a factored second-order truncation error that vanishes, which requires the extrapolation to be linear with quadratically-exact gradients. 

In this work, we incorporate recent advancements in the third-order edge-based scheme. First, we compute the lumped directed-area vector ${\bf n}_{jk}$ by an efficient algorithm proposed in Ref.~\cite{Nishikawa:RAA2025}, which does not require to the actual formation the medial dual control volume in a code. Second, the gradients are computed by the implicit edge-based algorithm described in Ref.~\cite{Nishikawa_aiaa_scitech2024}, which has been shown to be effective in stabilizing a nonlinear solver and also achieving lower errors in the numerical solution. Finally, we incorporate an accuracy-preserving boundary flux quadrature required for third-order accuracy \cite{nishikawa_boundary_quadrature:JCP2015} in an edge loop such that we only need to loop over boundary nodes (not elements) to close the residual, as described in Ref.~\cite{Nishikawa_aiaa_aviation2025_boundary}. 

Our focus is on the numerical flux. In the past, the third-order edge-based scheme has been used with a numerical flux in the form (\ref{upwind_flux}). Other numerical fluxes, which are not necessarily in the form (\ref{upwind_flux}), can be employed if they can be cast in the form (\ref{upwind_flux}). This is often possible, for example, as described for the HLLC flux function in Ref.~\cite{Fleischmann_LowDissCarbuncle_HLLC:JCP2020}. However, it is desirable to be able to directly use existing numerical flux functions. 

\vspace{-0.4cm}
\section{Flux Correction Form for a General Numerical Flux}
\label{intro}
\vspace{-0.2cm}

\indent

Before introducing a new approach, we recall the mechanism by which the edge-based scheme achieves third-order accuracy. Following the analysis in Appendix A of Ref.~\cite{nishikawa_liu_source_quadrature:jcp2017}, we consider the averaged flux part of the upwind flux (\ref{upwind_flux}) and expand it around node $j$ on an arbitrary tetrahedral grid, 
\vspace{-0.2cm}
\begin{eqnarray}
   \frac{1}{2} \left[  {\bf f}_L + {\bf f}_R   \right] 
   =
    {\bf f}_j + \frac{1}{2} \partial_{jk}  {\bf f}_j - \frac{1}{24} \partial_{jk}^3  {\bf f}_j  + O(h^4), 
    \quad
      \partial_{jk}  \equiv \Delta {\bf x}_{jk}  \cdot \left(  \partial_x, \partial_y, \partial_z \right),
      \label{ave_flux_expanded}
\end{eqnarray}
where $h$ is a typical mesh spacing and the derivatives are all evaluated at $j$. Then, expanding the source quadrature formula with $D=3$ by $ {\bf s}_k  =  {\bf s}_j + \partial_{jk}  {\bf s}_j  + \partial_{jk}^2  {\bf s}_j /2 + O(h^3)$, we find (see Ref.~\cite{nishikawa_liu_source_quadrature:jcp2017})
\begin{eqnarray}
\frac{1}{V_j}   {\bf Res}_j
  &=&
 \frac{1}{V_j}  \sum_{k \in \{ k_j\} }   \frac{1}{2} \left[     {\bf f}_L + {\bf f}_R   \right]    | {\bf n}_{jk} |  
  -  \frac{1}{60 V_j}  \sum_{k \in \{ k_j\} } \left[ 13   {\bf s}_j  + 3 ( \nabla  {\bf s}_j \cdot  \Delta {\bf x}_{jk} ) -  3 {\bf s}_k \right] \left( \Delta {\bf x}_{jk} \cdot {\bf n}_{jk}  \right) 
  \nonumber \\ [2ex]
    &=&
     \mbox{div}{\cal F}_j -  {\bf s}_j 
  - \frac{1}{24 V_j}  \sum_{k \in \{ k_j\} }   \partial_{jk}^3  {\bf f}_j    | {\bf n}_{jk} |   
  +  \frac{1}{40 V_j}     \sum_{k \in \{ k_j\} }   \partial_{jk}^2  {\bf s}_j  
\left( \Delta {\bf x}_{jk} \cdot {\bf n}_{jk}  \right)  + O(h^3) ,
\end{eqnarray}
where $V_j$ is the median dual volume. Note that the dissipation term has been ignored because it does not affect the leading second-order error \cite{nishikawa_liu_source_quadrature:jcp2017}. For third-order accuracy, there are two requirements: (1) the leading error is $O(h^2)$ on an arbitrary tetrahedral grid, and (2) the second-order error is factored, and thus, it vanishes on regular tetrahedral grids. The first requirement is met because the leading error in the above equation is $O(h^2)$. The second requirement is also met because the special source term quadrature formula employed here is designed for this purpose (see Ref.~\cite{nishikawa_liu_source_quadrature:jcp2017}). For example, for a regular tetrahedral grid constructed by dividing a hexahedral element into six tetrahedra, the above equation reduces to 
\vspace{-0.15cm}
\begin{eqnarray}
\frac{1}{V_j}   {\bf Res}_j
=
     \mbox{div}{\cal F}_j -  {\bf s}_j 
     - \frac{h^2}{12} 
     \left[ \partial_{xx} + \partial_{yy} + \partial_{zz} - \partial_{xy} - \partial_{yz} - \partial_{zx}  \right] (  \mbox{div}{\cal F}_j -  {\bf s}_j  ) + O(h^3).
\end{eqnarray}
If the smooth functions we expanded are the exact flux and source terms, then we have $\mbox{div}{\cal F}_j -  {\bf s}_j = 0$. Therefore, the first and second terms vanish, and we are left with a third-order truncation error. This completes the proof that the edge-based scheme achieves third-order accuracy on arbitrary tetrahedral grids. Readers are referred to Ref.~\cite{nishikawa_liu_source_quadrature:jcp2017} for further details. 


To enable the direct use of a general flux function, we consider the following form of the discretization:
\vspace{-0.15cm}
\begin{eqnarray}
  \sum_{k \in \{ k_j\} } \left[
{\bf f}_{jk}  +  \delta {\bf f}_{jk}   \right] | {\bf n}_{jk} |  =  \frac{1}{60}  \sum_{k \in \{ k_j\} } \left[ 13   {\bf s}_j  + 3  ( \nabla  {\bf s}_j \cdot \Delta {\bf x}_{jk} )   -  3 {\bf s}_k \right] \left( \Delta {\bf x}_{jk} \cdot {\bf n}_{jk}  \right) , 
\label{threed_fv_semidiscrete_system_00_FC}
\end{eqnarray}
where ${\bf f}_{jk}$ is the flux at the edge midpoint and $\delta {\bf f}_{jk}$ is a flux correction defined as
\begin{eqnarray}
  \delta {\bf f}_{jk}  =  \frac{C}{2} \left[   \left( \frac{\partial {\bf f} }{\partial {\bf w}} \right)_{\!\! j} \nabla {\bf w}_j 
  - \left( \frac{\partial {\bf f} }{\partial {\bf w}} \right)_{\!\! k} \nabla {\bf w}_k
   \right]  \cdot \Delta {\bf x}_{jk} ,
   \label{FC}
\end{eqnarray}
where the constant $C$ must be 1/4 to preserve third-order accuracy. To show this, we expand the flux ${\bf f}_{jk}$ and the flux correction $\delta {\bf f}_{jk}$ around node $j$:
\vspace{-0.25cm}
\begin{eqnarray}
 \!\!\!  \!\!\! 
{\bf f}_{jk}  \!\!  \!\!\! &=&   \!\!\! \!\!  {\bf f}_{j} + \frac{1}{2} \partial_{jk}  {\bf f}_j  + \frac{1}{8} \partial_{jk}^2  {\bf f}_j + \frac{1}{48} \partial_{jk}^3  {\bf f}_j + O(h^4),
 \\ [2ex]
  \!\!\!  \!\!\! 
  \delta {\bf f}_{jk}   \!\!\! \!\!  &=&  \!\!\!  \!\!
   \frac{C}{2} \left[  \nabla {\bf f}_j \!  - \!  \nabla {\bf f}_k 
   \right]  \! \cdot  \! \Delta {\bf x}_{jk} 
 \!  =\!
   \frac{C}{2} \! \left[   - \partial_{jk}  \left(  \nabla {\bf f}  \right)_j \!  - \frac{1}{2}  \partial_{jk}^2 \left(  \nabla {\bf f}  \right)_j \! + \! O(h^3) 
   \right] \!  \cdot \! \Delta {\bf x}_{jk} 
 \!  =\!
     \frac{C}{2} \! \left[   - \partial_{jk}^2 {\bf f}_j \! -  \! \frac{1}{2}  \partial_{jk}^3    \nabla {\bf f}_j  
   \right]   
\!+\! O(h^4)
   ,
\end{eqnarray}
and find
\begin{eqnarray}
{\bf f}_{jk}  +   \delta {\bf f}_{jk} 
=
{\bf f}_{j} + \frac{1}{2} \partial_{jk}  {\bf f}_j  +\left( \frac{1}{8} - \frac{C}{2}  \right) \partial_{jk}^2  {\bf f}_j + \left(  \frac{1}{48} -  \frac{C}{4}  \right) \partial_{jk}^3  {\bf f}_j 
+ O(h^4),
\end{eqnarray}
which matches Equation (\ref{ave_flux_expanded}) if we set $C=1/4$:
\begin{eqnarray}
{\bf f}_{jk}  +   \delta {\bf f}_{jk} 
=
{\bf f}_{j} + \frac{1}{2} \partial_{jk}  {\bf f}_j  + \left(  \frac{1}{48} -  \frac{1}{16}  \right) \partial_{jk}^3  {\bf f}_j  + O(h^4)
= 
{\bf f}_{j} + \frac{1}{2} \partial_{jk}  {\bf f}_j  -  \frac{1}{24}   \partial_{jk}^3  {\bf f}_j  + O(h^4).
\end{eqnarray}
Therefore, the same proof applies and thus, it preserves third-order accuracy. Finally, approximating ${\bf f}_{jk}$ by a general flux function ${\Phi}_{jk} ( {\bf u}_L, {\bf u}_R,  \hat{\bf n}_{jk} )$, we arrive at the following third-order edge-based discretization:
\begin{eqnarray}
  \sum_{k \in \{ k_j\} } \left[
{\Phi}_{jk} ( {\bf u}_L, {\bf u}_R,  \hat{\bf n}_{jk} ) +  \delta {\bf f}_{jk}   \right] | {\bf n}_{jk} |  =  \frac{1}{30}  \sum_{k \in \{ k_j\} } \left[ 13   {\bf s}_j  + 3  \nabla  {\bf s}_j \cdot  ( {\bf x}_k - {\bf x}_j )   -  3 {\bf s}_k \right] \left( \Delta {\bf x}_{jk} \cdot {\bf n}_{jk}  \right) , 
\label{threed_fv_semidiscrete_system_00_FC_practical}
\end{eqnarray}
where the flux correction is given by
\vspace{-0.2cm}
\begin{eqnarray}
  \delta {\bf f}_{jk}  =  \frac{1}{8} \left[   \left( \frac{\partial {\bf f} }{\partial {\bf w}} \right)_{\!\! j} \nabla {\bf w}_j 
  - \left( \frac{\partial {\bf f} }{\partial {\bf w}} \right)_{\!\! k} \nabla {\bf w}_k
   \right]  \cdot \Delta {\bf x}_{jk} .
   \label{FC}
\end{eqnarray}
It is emphasized that the numerical flux approximates the flux at the edge midpoint, and it must be evaluated with the left and right states quadratically extrapolated to the edge midpoint. This can be achieved without second derivatives by using the U-MUSCL scheme (\ref{umuscl_w}) with $\kappa=1/2$, which is exact as a solution extrapolation formula, for a quadratic function on an arbitrary grid, provided the gradient is exact for a quadratic function, as proved in Ref.~\cite{PadwayNishikawa:AIAAJ2022}. Therefore, the parameter $\kappa$ is not a free parameter in this case. The dissipation level, if desired, needs to be adjusted by other means, for example, by the choice of the numerical flux function or by extending a technique used to develop a low-dissipation numerical flux for the NASA FUN3D code, as described in Ref.~\cite{nishikawa_liu_aiaa2018-4166}.

\vspace{-0.45cm}
\section{Results}
\label{results}
\vspace{-0.2cm}

\indent

\begin{figure}[t] 
    \centering
  \begin{subfigure}[t]{0.49\textwidth}
  \includegraphics[width=0.6\textwidth,trim=2 2 2 2 ,clip]{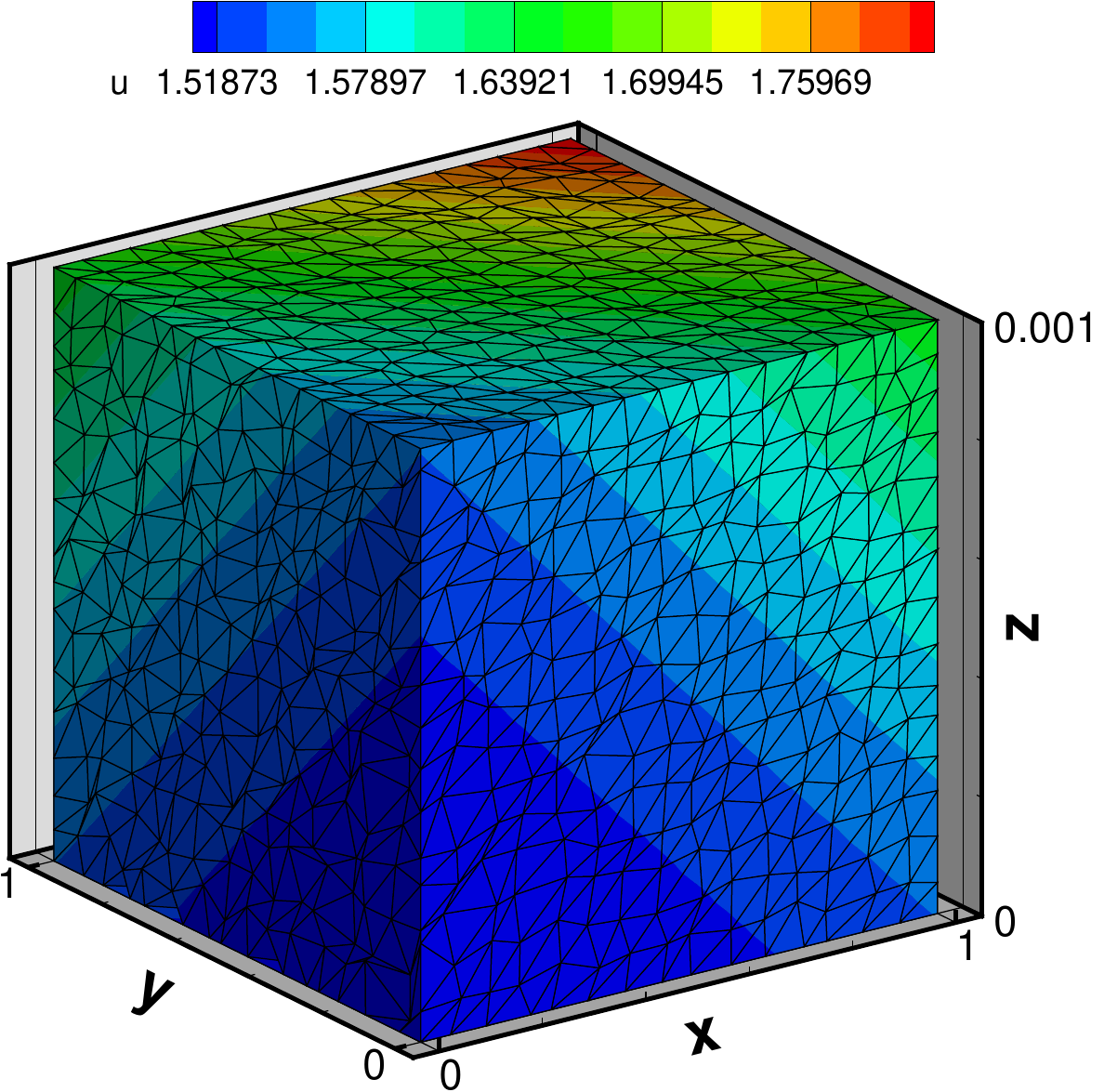}
  \caption[]{Coarsest grid with $x$-velocity contours. 
  \label{fig:accuracy_verification_grid}} 
  \end{subfigure}
  \begin{subfigure}[t]{0.49\textwidth}
  \includegraphics[width=0.6\textwidth,trim=0 0 0 0 ,clip]{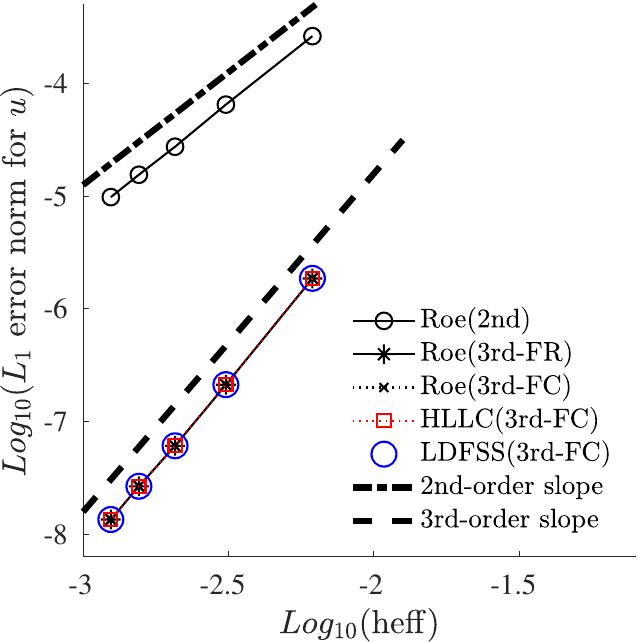}
  \caption[]{Error convergence results. 
  \label{fig:accuracy_verification_err}
  } 
  \end{subfigure} 
\caption{
Accuracy verification test: the coarsest grid and error convergence.
 \label{fig:accuracy_verification}}
\end{figure} 

For the demonstration that the proposed flux-correction form preserves third-order accuracy, it suffices to perform an accuracy verification study using the method of 
manufactured solutions. Thus, we consider the following vector of functions: ${\bf w} = {\bf w}_0 + \delta {\bf w} \exp( 0.2 x + 0.2 y + 200 z  )$, where ${\bf w}_0 = (1.0, 0.3, 0.2, 0.1, 1.0)$, $\delta {\bf w} = (0.1, 0.1, 0.1, 0.1, 0.1)$, and then compute the forcing terms ${\bf s}({\bf x})$ by substituting them in the Euler equations, so that the above functions serve as exact solutions. We solve the residual system by an implicit defect-correction solver for a series of consistently refined irregular tetrahedral grids with $n$$\times$$n$$\times$$n$ nodes, where $n=16, 32, 48, 64$, and $80$ in a high-aspect-ratio domain $(x,y,z) \in [0,1] \times [0,1] \times [0,0.001]$. The coarsest mesh and solution contours are shown in Figure~\ref{fig:accuracy_verification_grid}, where the $z$-coordinate is scaled by 1000, so that the domain looks like a cube. In reality, the domain is thin and the grid is highly skewed. The solver is taken to be converged when the $L_1$ norms of the residuals 
are reduced by six orders of magnitude from the initial value evaluated with randomly-perturbed exact solutions specified at all nodes. Discretization errors are measured in the $L_1$ norm against exact point-valued solutions at nodes. Boundary conditions are imposed weakly through a numerical flux with the right state specified by exact solution values (see Ref.~\cite{Nishikawa_aiaa_aviation2025_boundary} for details).

 As examples of a general flux function not in the form (\ref{upwind_flux}), we consider the HLLC flux \cite{HLLC_Batten_SIAM1997} and the LDFSS flux \cite{Edwards_LDFSS_1997} implemented in the proposed flux-correction form (\ref{threed_fv_semidiscrete_system_00_FC_practical}). These fluxes will be referred to as HLLC(3rd-FC) and LDFSS(3rd-FC). Also, we consider the Roe flux  \cite{Roe_JCP_1981} implemented in the original form (\ref{upwind_flux}) and the proposed flux-correction form (\ref{threed_fv_semidiscrete_system_00_FC_practical}). The former will be referred to as Roe(3rd-FR) and the latter will be referred to as Roe(3rd-FC). Finally, we also consider a second-order accurate scheme with the Roe flux, which will be referred to as Roe(2nd). For the second-order scheme, we use $\kappa=1/2$ for the solution extrapolation, and the fluxes are evaluated with the extrapolated solutions: ${\bf f}_L = {\bf f}({\bf w}_L)$ and ${\bf f}_R = {\bf f}({\bf w}_R)$. The source term is simply evaluated at a node, not by the accuracy-preserving source quadrature formula.

Error convergence results are shown for the $x$-velocity component in Figure \ref{fig:accuracy_verification_err} (results are similar for other variables, and thus not shown). As can be seen, the third-order edge-based scheme achieves third-order accuracy successfully with all the numerical flux functions, and it is significantly more accurate than the second-order scheme.

\vspace{-0.35cm}
\section{Conclusions}
\label{conclusions}
\vspace{-0.2cm}

\indent

In this short note, we have proposed a flux-correction form of the third-order edge-based scheme, as an alternative to flux extrapolation. The proposed method enables the direct use of a general numerical flux function and is also expected to simplify the implementation of the third-order edge-based scheme in an existing edge-based solver, at least in the numerical flux part, because one only needs to use the U-MUSCL scheme with $\kappa=1/2$ and add a flux-correction term. Numerical results show that third-order accuracy is preserved with the proposed method using the HLLC and LDFSS flux functions in their original forms. The impact on practical applications will be investigated in the future using practical edge-based solvers, and a more meaningful comparative study with a second-order edge-based solver can be performed because because the same numerical flux can be used for both second- and third-order edge-based solvers.

\vspace{-0.35cm}
\section*{Acknowledgments}
\vspace{-0.2cm}
 

The author gratefully acknowledges support by NASA Langley Research Center under Contract No.80LARC23DA003.


\vspace{-0.35cm}
\bibliography{../../bibtex_nishikawa_database}
\bibliographystyle{unsrt}

\end{document}